\documentclass[11pt]{amsart}
\usepackage[all]{xy}
\usepackage{amsmath}
\usepackage{amsfonts}
\usepackage{amssymb}
\usepackage{latexsym}
\usepackage{graphicx}
\usepackage{color}
\usepackage{amscd}
\usepackage{epsfig}
\usepackage{cite}

\newtheorem{problem}{Problem}[section]
\newtheorem{definition}[problem]{Definition}
\newtheorem{lemma}[problem]{Lemma}
\newtheorem{theorem}[problem]{Theorem}

\newtheorem{conjecture}[problem]{Conjecture}

\title{Generic $T$-adic exponential sums in one variable}
\author{Chunlei Liu}
\address{Department of Mathematics, Shanghai Jiao Tong
University, Shanghai 200240, P.R. China, E-mail: clliu@sjtu.edu.cn}
\author{Wenxin Liu}\address{School of Mathematical Sciences, Beijing Normal
University, Beijing 100875, P.R. China, E-mail:
wenxin8210@mail.bnu.edu.cn}\author{Chuanze Niu}\address{School of
Mathematical Sciences, Beijing Normal University, Beijing 100875,
P.R. China, E-mail: czniubnu@yahoo.cn}
\thanks{This research is supported by NSFC Grant No.
10671015.}
\begin{document}
\maketitle
\begin{abstract}
The $T$-adic exponential sum associated to a Laurent polynomial in
one variable is studied. An explicit arithmetic polygon is proved to
be the generic Newton polygon of the $C$-function of the T-adic
exponential sum. It gives the generic Newton polygon of
$L$-functions of $p$-power order exponential sums.
\end{abstract}



\section{Introduction}

Let $W$ be the ring scheme of Witt vectors, $\mathbb{F}_q$ the field
of characteristic $p$ with $q$ elements,
$\mathbb{Z}_{q}=W(\mathbb{F}_q)$, and
$\mathbb{Q}_q=\mathbb{Z}_q[\frac{1}{p}]$.

Let $\triangle\supsetneq\{0\}$ be an integral convex polytope
 in $\mathbb{R}^n$, and $I$ the set of vertices of
 $\triangle$ different from the origin.
Let
$$f(x)=\sum\limits_{u\in \triangle}(a_ux^u,0,0,\cdots)\in
W(\mathbb{F}_q[x_1^{\pm1},\cdots,x_n^{\pm1}])\text{ with }
\prod_{u\in I}a_u\neq0,$$ where $x^u=x_1^{u_1}\cdots x_n^{u_n}$ if
$u=(u_1,\cdots,u_n)\in\mathbb{Z}^n$.

Let $T$ be a variable.\begin{definition} The sum
$$S_{f}(k,T)=\sum\limits_{x\in(\mathbb{F}_{q^k}^{\times})^n}
(1+T)^{{\rm
Tr}_{\mathbb{Z}_{q^k}/\mathbb{Z}_p}(f(x))}\in\mathbb{Z}_p[[T]]$$ is
call a $T$-adic exponential sum. And the function
$$L_f(s,T)=\exp(\sum\limits_{k=1}^{\infty}S_f(k,T)\frac{s^k}{k})\in
1+s\mathbb{Z}_p[[T]][[s]]$$ is called an $L$-function of $T$-adic
exponential sums. \end{definition}We view $L_f(s, T)$ as a power
series in the single variable $s$ with coefficients in the $T$-adic
complete field $\mathbb{Q}_p((T))$. Let $\zeta_{p^m}$ be a primitive
$p^m$-th root of unity, and $\pi_m=\zeta_{p^m}-1$. Then
$L_f(s,\pi_m)$ is the $L$-function of the $p$-power order
exponential sums $S_f(k,\pi_m)$ studied by Liu-Wei \cite{LW}.

\begin{definition}The function
$$C_f(s,T) =C_f(s,T;\mathbb{F}_q) =\exp(\sum\limits_{k=1}^{\infty}-(q^k-1)^{-n}S_{f}(k,T)\frac{s^k}{k})$$
is called a $C$-function of $T$-adic exponential
sums.\end{definition} We have
$$L_f(s,T) = \prod_{i=0}^n C_f(q^is, T)^{(-1)^{n-i+1}{n\choose i}},$$
and $$C_f(s, T)= \prod_{j=0}^{\infty} L_f(q^js,
T)^{(-1)^{n-1}{n+j-1\choose j}}.$$ So the $C$-function $C_f(s,T)$
and the $L$-function $L_f(s,T)$ determine each other. From the last
identity, one sees that
$$C_f(s,T)\in
1+s\mathbb{Z}_p[[T]][[s]].$$We also view $C_f(s, T)$ as a power
series in the single variable $s$ with coefficients in the $T$-adic
complete field $\mathbb{Q}_p((T))$. The $C$-function $C_f(s,T)$ was
shown to be $T$-adic entire in $s$ by Liu-Wan \cite{LWn}.

Let $C(\triangle)$ be the cone generated by $\triangle$, and
$M(\triangle)=M(\triangle)\cap \mathbb{Z}^n$. There is a degree
function $\deg$ on $C(\triangle)$ which is
$\mathbb{R}_{\geq0}$-linear and takes the values $1$ on each
co-dimension $1$ face not containing $0$. For $a\not\in
C(\triangle)$, we define $\deg(a)=+\infty$.
\begin{definition}A convex function on
$[0,+\infty]$ which is linear between consecutive integers with
initial value $0$ is called the infinite Hodge polygon of
$\triangle$ if its slopes between consecutive integers are the
numbers $\deg(a)$, $a\in M(\triangle)$. We denote this polygon by
$H_{\triangle}^{\infty}$.
\end{definition}

Liu-Wan \cite{LWn} also proved the following.

\begin{lemma}We have $$T-\text{adic NP of }C_{f}(s,T)\geq\text{ord}_p(q)(p-1)H_{\triangle}^{\infty},$$
where NP stands for Newton polygon.\end{lemma}

\begin{definition}The $T^a$-adic Newton polygon of $C_f(s,T;\mathbb{F}_{p^a})$ is called
the absolute Newton polygon of
$C_f(s,T;\mathbb{F}_{p^a})$.\end{definition}
\begin{conjecture}\label{gent}If $p$ is sufficiently large, then the
absolute $T$-adic Newton polygon of $C_f(s,T)$ is constant for a
generic $f$. We call it the generic Newton polygon of $C_f(s,T)$.
\end{conjecture}
\begin{definition}The $\pi_m^a$-adic Newton polygon of $C_f(s,\pi_m;\mathbb{F}_{p^a})$ is called
the absolute Newton polygon of
$C_f(s,\pi_m;\mathbb{F}_{p^a})$.\end{definition} Combine results of
Gelfand-Kapranov-Zelevinsky \cite{GKZ}, Adolphson-Sperber \cite{AS},
Liu-Wei \cite{LW} with Grothendieck specialization lemma \cite{Ka},
the absolute Newton polygon of $C_f(s,\pi_m)$ is constant for a
generic $f$. We call it the generic Newton polygon of
$C_f(s,\pi_m)$.
\begin{conjecture}\label{gen} If $p$ is sufficiently large, then
the generic Newton polygon of $C_f(s,\pi_m)$ is independent of $m$,
and coincides with the generic Newton polygon of $C_f(s,T)$.
\end{conjecture}
The generic Newton polygon of $C_f(s,\pi_m)$ for $m=1$ was studied
by Wan \cite{Wa1, Wa2}.

In the rest of this section we assume that
$\triangle\subset\mathbb{Z}$.
\begin{definition}Let $0\neq a\in M(\triangle)$. We define
$$\delta_{\in}(a)=\left\{
                           \begin{array}{ll}
                            1, & \hbox{ }\{\deg(a)\}=\{\deg(pi)\}\text{ for some }i\text{ with }ia>0,\deg(i)<\{\deg(a)\}, \\
                            0, & \hbox{ }{\rm otherwise},
                            \end{array}
                            \right.
$$
where $\{\cdot\}$ is the fractional part of a real number. We also
define $\delta_{\in}(0)=0$.
\end{definition}

\begin{definition}A convex function on
$[0,+\infty]$ which is linear between consecutive integers with
initial value $0$ is called the arithmetic polygon of $\triangle$ if
its slopes between consecutive integers are the numbers
$$\varpi_{\triangle}(a)=\lceil(p-1)\deg(a)\rceil-\delta_{\in}(a), a\in M(\triangle),$$
where $\lceil\cdot\rceil$ is the least integer equal or greater than
a real number. We denote this polygon by $p_{\triangle}$.
\end{definition}

We can prove the following.
\begin{theorem}\label{arithhodge}We
have$$p_{\triangle}\geq (p-1)H_{\triangle}^{\infty}.$$ Moreover,
they coincide at the point ${\rm Vol}(\triangle)$.
\end{theorem}

Let $D$ be the least common multiple of the nonzero endpoint(s) of
$\triangle$. The main results of this paper are the following
theorems.
\begin{theorem}\label{main1}If $p>3D$, then
$$T-\text{adic NP of }C_{f}(s,T)\geq {\rm ord}_p(q)p_{\triangle}.$$
\end{theorem}
\begin{theorem}\label{main2}Let $f(x)=\sum\limits_{u\in\triangle}(a_ux^u,0,0,\cdots)$, and $p>3D$.
Then there is a non-zero polynomial $H(y)\in\mathbb{F}_q[y_u\mid
u\in\triangle]$ such that $$T-\text{adic NP of }C_{f}(s,T)={\rm
ord}_p(q)p_{\triangle}$$ if and only if
$H((a_u)_{u\in\triangle})\neq0$.
\end{theorem} The above theorem implies Conjecture
\ref{gent} for $\triangle\subset\mathbb{Z}$.
\begin{theorem}\label{main3}Let $f(x)=\sum\limits_{u\in\triangle}(a_ux^u,0,0,\cdots)$, $p>3D$, and $m\geq1$.
Then $$\pi_m-\text{adic NP of }C_{f}(s,\pi_m)={\rm
ord}_p(q)p_{\triangle}$$ if and only if
$H((a_u)_{u\in\triangle})\neq0$.
\end{theorem}
 The above  theorem implies Conjecture
\ref{gen} for $\triangle\subset\mathbb{Z}$.

Note that, if $p\nmid D$, then $L(s,\pi_m)$ is a polynomial of
degree $p^{m-1}{\rm Vol}(\triangle)$ for all $m\geq1$ by a result of
Adolphson-Sperber \cite{AS} and a result of Liu-Wei \cite{LW}. We
shall prove the following.
\begin{theorem}\label{expsum}Let $f(x)=\sum\limits_{u\in\triangle}(a_ux^u,0,0,\cdots)$, and $p>3D$.
Then $$\pi_m-\text{adic NP of }L_{f}(s,\pi_m)\geq{\rm
ord}_p(q)p_{\triangle}\text{ on }[0,p^{m-1}{\rm Vol}(\triangle)]$$
with equality holding if and only if
$H((a_u)_{u\in\triangle})\neq0$.
\end{theorem}
Zhu \cite{Zh1, Zh2} and Blache-F\'{e}rard \cite{BF} studied the
Newton polygon of the $L$-function $L_{f}(s,\pi_m)$  for $m=1$.
\section{Arithmetic estimate}
In this section $\triangle\subset\mathbb{Z}$, $A$ is a finite subset
of $M(\triangle)\times\mathbb{Z}/(b),$ and $\tau$ is a permutation
of $A$. We shall estimate
$$\sum\limits_{a\in A}\lceil
\deg(pa-\tau(a)\rceil,$$ where $\deg(i,u)=\deg(i)$.

However, except in the ending paragraph, we assume that
$\triangle=[0,d]$, $A$ is a finite subset of
$\{1,2,\cdots\}\times\mathbb{Z}/(b)$.

Write
$$\{x\}'=1+x-\lceil x\rceil=\left\{
                                                   \begin{array}{ll}
                                                     \{x\}, & \hbox{ if } \{x\}\neq0,\\
                                                     1, & \hbox{ if
}\{x\}=0.
                                                   \end{array}
                                                 \right.$$
\begin{lemma}We have$$\sum\limits_{a=0}^{m}(\delta_{<}-\delta_{\in})(a)=\#\{1\leq a\leq
d\{\frac{m}{d}\}'\mid \{\frac{m}{d}\}'<\{\frac{pa}{d}\}'\},$$ where
$$\delta_{<}(a)=\left\{
                  \begin{array}{ll}
                    1, & \hbox{ } \{\frac{a}{d}\}'<\{\frac{pa}{d}\}'\\
                    0, & \hbox{ }{\rm otherwise}.
                  \end{array}
                \right.$$\end{lemma}\proof Note that
both $\delta_{\in}$ and $\delta_{<}$ have a period $d$ and have
initial value $0$. So we may assume that $m<d$. We have
$$\sum\limits_{a=0}^m\delta_{\in}(a)=\sum\limits_{a=1}^m
\sum\limits_{\stackrel{i=1}{pi\equiv a(
d)}}^{a-1}1$$$$=\sum\limits_{i=1}^{m-1}
\sum\limits_{\stackrel{a=i+1}{pi\equiv a( d)}}^{m}1=\#\{1\leq
i<m\mid i<d\{\frac{pi}{d}\}\leq m\}.$$And, by definition,
$$\sum\limits_{a=0}^m\delta_{<}(a)=\#\{1\leq a\leq m\mid
a<d\{\frac{pa}{d}\}\}.$$The lemma now follows.
\endproof
\begin{lemma}Let $A,B,C$ be sets
with $A$ finite, and $\tau$ a permutation of $A$. Then
$$\#\{a\in A\mid \tau(a)\in B,a\in C\}
$$$$\geq\#\{a\in A\mid a\in B,a\in C\}-\#\{a\in A\mid a\not\in
B,a\not\in C\}.$$
\end{lemma}
\proof We have$$\#\{a\in A \mid a\in B,\tau(a)\in B,a\in C\}
$$$$\geq\#\{a\in A \mid a\in B,a\in C\}-\#\{a\in A\mid a\in
B,\tau(a)\not\in B\},$$ and
$$\#\{a\in A\mid a\not\in B,\tau(a)\in B,a\in C\} $$$$\geq\#\{a\in
A\mid a\not\in B,\tau(a)\in B\}-\#\{a\in A\mid a\not\in B,a\not\in
C\}.$$ Note that
$$\#\{a\in A\cap B\mid \tau(a)\not\in B\}=\#\{a\in A\setminus B\mid \tau(a)\in B\}.$$ So
$$\#\{a\in A \mid \tau(a)\in B,a\in C\}
$$$$\geq\#\{a\in A\mid a\in B,a\in C\}-\#\{a\in A\mid a\not\in
B,a\not\in C\}.$$The lemma is proved.\endproof For $a\in A$, we
define
$$\delta_{<}^{\tau}(a)=\left\{
                         \begin{array}{ll}
                           1, & \hbox{ } \{\deg(\tau(a))\}'<\{p\deg(a)\}',\\
                           0, & \hbox{otherwise.}
                         \end{array}
                       \right.$$
\begin{theorem}If $p>3d$, and
$$A=(\{1,\cdots,m-1\}\times\mathbb{Z}/(b))\cup\{(m,i_0),\cdots,(m,i_{l-1})\},$$
then$$\sum\limits_{a\in A}\delta_{<}^{\tau}(a)\geq
b\sum\limits_{a=0}^{m-1}(\delta_{<}-\delta_{\in})(a)+l(\delta_{<}-\delta_{\in})(m).$$
\end{theorem}\proof
First we assume that $l=0$. We have
$$\sum\limits_{a\in
A}\delta_{<}^{\tau}(a)\geq\#\{a\in
A,\{\deg(\tau(a))\}'\leq\{\frac{m-1}{d}\}'<\{p\deg(a)\}'\}.$$
Applying the last lemma with
$$B=\{a\in
A,\{\deg(a)\}'\leq\{\frac{m-1}{d}\}'\},$$ and $$ C=\{a\in
A,\{\frac{m-1}{d}\}'<\{p\deg(a)\}'\},$$ we get
$$\sum\limits_{a\in
A}\delta_{<}^{\tau}(a)\geq\#\{a\in
A,\{\deg(a)\}'\leq\{\frac{m-1}{d}\}'<\{p\deg(a)\}'\}$$$$-\#\{a\in
A,\{\deg(a)\}'>\{\frac{m-1}{d}\}'\geq\{p\deg(a)\}'\}.$$ We have
$$\#\{a\in
A,\{\deg(a)\}'\leq\{\frac{m-1}{d}\}'<\{p\deg(a)\}'\}$$
$$=b([\frac{m}{d}]+1)\#\{1\leq a\leq d\{\frac{m-1}{d}\}'\mid \{\frac{m-1}{d}\}'<\{\frac{pa}{d}\}'\}.$$
We also have
$$\#\{a\in
A,\{\deg(a)\}'>\{\frac{m-1}{d}\}'\geq\{p\deg(a)\}'\}$$
$$=b[\frac{m}{d}]\#\{d\geq a>d\{\frac{m-1}{d}\}'\mid \{\frac{m-1}{d}\}'\geq\{\frac{pa}{d}\}'\}.$$
$$=b[\frac{m}{d}]\#\{1\leq a\leq d\{\frac{m-1}{d}\}'\mid \{\frac{m-1}{d}\}'<\{\frac{pa}{d}\}'\}.$$

It follows that
$$\sum\limits_{a\in
A}\delta_{<}^{\tau}(a)\geq b\#\{1\leq a\leq d\{\frac{m-1}{d}\}'\mid
\{\frac{m-1}{d}\}'<\{\frac{pa}{d}\}'\} $$$$\geq
b\sum\limits_{a=0}^{m-1}(\delta_{<}-\delta_{\in})(a).$$

Secondly we assume that $\delta_{\in}(m)=0$. Extend  the action of
$\tau$ trivially from $A$ to $\{1,\cdots,m\}\times\mathbb{Z}/(b)$.
By what we just proved, $$\sum\limits_{a\in
\{1,\cdots,m\}\times\mathbb{Z}/(b)}\delta_{<}^{\tau}(a)\geq
b\sum\limits_{a=0}^{m}(\delta_{<}-\delta_{\in})(a).$$ It follows
that$$\sum\limits_{a\in A}\delta_{<}^{\tau}(a)\geq
b\sum\limits_{a=0}^{m}(\delta_{<}-\delta_{\in})(a)-(b-l)\delta_{<}(m)$$$$\geq
b\sum\limits_{a=0}^{m-1}(\delta_{<}-\delta_{\in})(a)+l(\delta_{<}-\delta_{\in})(m).$$

Finally we assume that $\delta_{\in}(m)=1$. We have
$$\sum\limits_{a\in
A}\delta_{<}^{\tau}(a)\geq\#\{a\in
A,\{\deg(\tau(a))\}'\leq\{\frac{m-1}{d}\}'<\{p\deg(a)\}'\}.$$
Applying the last lemma with
$$B=\{a\in
A,\{\deg(a)\}'\leq\{\frac{m-1}{d}\}'\},$$ and $$ C=\{a\in
A,\{\frac{m-1}{d}\}'<\{p\deg(a)\}'\},$$ we get
$$\sum\limits_{a\in
A}\delta_{<}^{\tau}(a)\geq\#\{a\in
A,\{\deg(a)\}'\leq\{\frac{m-1}{d}\}'<\{p\deg(a)\}'\}$$$$-\#\{a\in
A,\{\deg(a)\}'>\{\frac{m-1}{d}\}'\geq\{p\deg(a)\}'\}.$$ We have
$$\#\{a\in
A,\{\deg(a)\}'\leq\{\frac{m-1}{d}\}'<\{p\deg(a)\}'\}$$
$$=b([\frac{m}{d}]+1)\#\{1\leq a\leq d\{\frac{m-1}{d}\}'\mid \{\frac{m-1}{d}\}'<\{\frac{pa}{d}\}'\}.$$
We also have
$$\#\{a\in
A,\{\deg(a)\}'>\{\frac{m-1}{d}\}'\geq\{p\deg(a)\}'\}$$
$$=b[\frac{m}{d}]\#\{d\geq a>d\{\frac{m-1}{d}\}'\mid \{\frac{m-1}{d}\}'\geq\{\frac{pa}{d}\}'\}
+l1_{\{\frac{m-1}{d}\}'\geq\{\frac{pm}{d}\}'}.$$
$$=b[\frac{m}{d}]\#\{1\leq a\leq d\{\frac{m-1}{d}\}'\mid \{\frac{m-1}{d}\}'<\{\frac{pa}{d}\}'\}
+l1_{\{\frac{m-1}{d}\}'\geq\{\frac{pm}{d}\}'}.$$

It follows that
$$\sum\limits_{a\in
A}\delta_{<}^{\tau}(a)\geq b\#\{1\leq a\leq d\{\frac{m-1}{d}\}'\mid
\{\frac{m-1}{d}\}'<\{\frac{pa}{d}\}'\}$$
$$-l1_{\{\frac{m-1}{d}\}'\geq\{\frac{pm}{d}\}'}\geq
b\sum\limits_{a=0}^{m-1}(\delta_{<}-\delta_{\in})(a)+l(\delta_{<}-\delta_{\in})(m).$$
The proof of the theorem is completed.\endproof

\begin{theorem}If $p>3d$, $A$ is of
cardinality $bm+l$ with $0\leq l<b$, then
$$\sum\limits_{a\in A}(\lceil
p\deg(a)\rceil-\lceil\deg(a)\rceil+\delta_{<}^{\tau}(a))\geq
bp_{\triangle}(m)+l\varpi(m).$$ Moreover, the strict inequality
holds if $A$ is not of the form
$$(\{1,\cdots,m-1\}\times\mathbb{Z}/(b))\cup\{(m,i_0),\cdots,(m,i_{l-1})\}.$$
\end{theorem}
\proof  We may assume that $A$ is not of the form
$$(\{1,\cdots,m-1\}\times\mathbb{Z}/(b))\cup\{(m,i_0),\cdots,(m,i_{l-1})\}.$$
There is an element $a_0\in A$ with $\deg(a_0)>\deg(m)$. Set
$A'=A\setminus\{a_0\}$ and
$$\tau'(a)=\left\{
             \begin{array}{ll}
               \tau(a), & \hbox{ } a\neq\tau^{-1}(a_0),\\
               \tau(a_0), & \hbox{}a=\tau^{-1}(a_0).
             \end{array}
           \right.
$$
Then
$$\sum\limits_{a\in A}(\lceil
p\deg(a)\rceil-\lceil\deg(a)\rceil+\delta_{<}^{\tau}(a))$$$$>\sum\limits_{a\in
A'}(\lceil
p\deg(a)\rceil-\lceil\deg(a)\rceil+\delta_{<}^{\tau'}(a))+\varpi(m).$$The
theorem now follows by induction.\endproof

We now assume that $\triangle\subset\mathbb{Z}$, and $A$ is a finite
subset of $M(\triangle)\times\mathbb{Z}/(b).$ For an integer
$m\geq1$, we define
$$A_m=\{a\in M(\triangle)\mid \varpi(a)\leq p_{\triangle}(m)-p_{\triangle}(m-1)\}.$$
\begin{theorem}\label{order-estimate}If $p>3D$, then
$$\sum\limits_{a\in A}\lceil\deg(pa-\tau(a))\rceil\geq bp_{\triangle}(m).$$
Moreover, the strict inequality holds if $m$ is a turning point of
$p_{\triangle}$, and $A\neq A_m\times\mathbb{Z}/(b)$.
\end{theorem}

\proof We define ${\rm sgn}((i,u))={\rm sgn}(i)$. Write
$$A_i=\{a\in A\mid {\rm sgn}(a)=(-1)^i\},$$ and $$A_{ij}=\{a\in A\mid {\rm sgn}(a)=(-1)^i,{\rm
sgn}(\tau(a))=(-1)^j\}.$$ Define a new permutation $\tau_0$ as
follows:\begin{itemize}
          \item $\tau_0$ is identity on $A_0$.
          \item $\tau_0=\tau$ on $A_{11}$ and $A_{22}$.
          \item $\tau_0$ maps $A_1\setminus A_{11}$ to
$A_1\setminus\tau(A_1)$.\item $\tau_0$ maps $A_2\setminus A_{22}$ to
$A_2\setminus\tau(A_2)$.
        \end{itemize}We have
$$\sum\limits_{a\in A}\lceil\deg(pa-\tau(a))\rceil\geq
\sum\limits_{a\in A}(\lceil
p\deg(a)\rceil-\lceil\deg(\tau_0(a))\rceil+1_{\{\deg(\tau_0(a))\}'<\{p\deg(a)\}'}).$$The
theorem now follows the last one.\endproof
\section{The $T$-adic Dwork Theory}
In this section we review the $T$-adic analogue of Dwork theory on
exponential sums.

Let
$$E(t)=\exp(\sum_{i=0}^{\infty}\frac{t^{p^i}}{p^i})=\sum\limits_{i=0}^{+\infty}\lambda_it^i
\in 1+t{\mathbb Z}_p[[t]]$$ be the $p$-adic Artin-Hasse exponential
series. Define a new $T$-adic uniformizer $\pi$ of ${\mathbb
Q}_p((T))$ by the formula $E(\pi)=1+T$. Let $\pi^{1/D}$ be a fixed
$D$-th root of $\pi$. Let
$$L=\{\sum_{i\in M(\triangle)}c_i\pi^{\deg(i)}x^i:\
 c_i\in\mathbb{Z}_q[[\pi^{1/D}]] \}.$$

Let $a\mapsto\hat{a}$  be the Teichm\"{u}ller lifting. One can show
that the series
$$E_f(x) :=\prod\limits_{a_i\neq0}E(\pi \hat{a}_ix^i)\in L.$$
Note that the Galois group of $\mathbb{Q}_q$ over
$\mathbb{Q}_p$ can act on $L$ but keeping $\pi^{1/D}$ as well as the
variable $x$ fixed.
 Let $\sigma$ be the Frobenius element in the
Galois group such that $\sigma(\zeta)=\zeta^p$ if $\zeta$ is a
$(q-1)$-th root of unity. Let $\Psi_p$ be the operator on $L$
defined by the formula
$$\Psi_p(
 \sum\limits_{i\in
M(\triangle)} c_ix^i)=\sum\limits_{i\in M(\triangle)} c_{pi}x^i.$$
Then $\Psi:=\sigma^{-1}\circ\Psi_p\circ E_f$ acts on the $T$-adic
Banach module
$$B=\{\sum\limits_{i\in M(\triangle)}c_i\pi^{\deg(i)}x^i \in L,\
 \text{\rm ord}_T(c_i)\rightarrow+\infty
 \text{ if }\deg(i)\rightarrow+\infty\}.$$
We call it Dwork's $T$-adic semi-linear operator because it is
semi-linear over $\mathbb{Z}_q[[\pi^{\frac{1}{D}}]]$ .

Let $b=\log_pq$. Then the $b$-iterate $\Psi^b$ is linear over
$\mathbb{Z}_q[[\pi^{1/D}]]$, since
$$\Psi^{b}=\Psi_p^{b}\circ
\prod\limits_{i=0}^{b-1}E_{f}^{\sigma^i}(x^{p^i}).$$ One can show
that $\Psi$ is completely continuous in the sense of Serre
\cite{Se}. So $\det(1-\Psi^bs\mid
B/\mathbb{Z}_{q}[[\pi^{\frac{1}{D}}]])$ and $\det(1-\Psi s\mid
B/\mathbb{Z}_p[[\pi^{\frac{1}{D}}]])$ are well-defined.

We now state the $T$-adic Dwork trace formula\cite{LWn}.
\begin{theorem}
We have
$$C_f(s,T)=\det(1-\Psi^bs\mid
B/\mathbb{Z}_{q}[[\pi^{\frac{1}{D}}]]).$$
\end{theorem}
\section{The Dwork semi-linear operator}

Write $$E_f(x)= \sum\limits_{i\in M(\triangle)}\gamma_ix^i,$$ and
$$\det(1-\Psi s\mid B/\mathbb{Z}_p[[\pi^{\frac{1}{D}}]])=\sum\limits_{i=0}^{+\infty}(-1)^ic_is^i.$$
 Let $O(\pi^{\alpha})$ denotes any element of $\pi$-adic order $\geq
\alpha$. In this section we prove the following.
\begin{theorem}\label{general-points}Let $p>3D$. Then
$$\text{ord}_{\pi}(c_{bm})\geq
bp_{\triangle}(m).$$ Moreover, if $m<{\rm Vol}(\triangle)$ is a
turning point of $p_{\triangle}$, then
$$c_{bm}=\pm{\rm Norm}(\det(\gamma_{pi-j})_{i,j\in A_m})+O(\pi^{bp_{\triangle}(m)+1/D}),$$
where Norm is the norm map from  $\mathbb{Q}_q(\pi^{1/D})$ to
$\mathbb{Q}_p(\pi^{1/D})$.
\end{theorem}

\proof Fix a normal basis $\bar{\xi}_u$, $u\in\mathbb{Z}/(b)$ of
$\mathbb{F}_q$ over $\mathbb{F}_p$. Let $\xi_u$ be their
Terchm\"{u}ller lift of $\bar{\xi}_u$. Then $\xi_u$,
$u\in\mathbb{Z}/(b)$ is a normal basis of $\mathbb{Q}_q$ over
$\mathbb{Q}_p$, and $\sigma$ acts on the basis $\xi_u$ ,
$u\in\mathbb{Z}/(b)$ as a permutation. Let
$(\gamma_{(i,u),(j,\omega)})_{i,j\in M(\triangle),1\leq u,\omega\leq
b}$ be the matrix of $\Psi$ on
$B\otimes_{\mathbb{Z}_p}\mathbb{Q}_p(\pi^{1/D})$ with respect to the
basis $\{\xi_ux^i\}_{i\in M(\triangle),1\leq u\leq b}$. Then
$$c_{bm}=\sum\limits_{A}\text{det}((\gamma_{(i,u),(j,\omega)})_{(i,u),(j,\omega)\in
A}),
$$
where $A$ runs over all subsets of $
M(\triangle)\times\mathbb{Z}/(b)$ with cardinality $bm$. One can
show that $$\det(\gamma_{i,j})_{i,j\in
A_m\times\mathbb{Z}/(b)}=\pm{\rm Norm}(\det((\gamma_{pi-j})_{i,j\in
A_m})).$$Therefore the theorem follows from the following.\endproof
\begin{theorem}\label{det}
Let $A\subset M(\triangle)\times\mathbb{Z}/(b)$ be a subset of
cardinality $bm$. If $p>3D$, then
$${\rm ord}_T(\det(\gamma_{(i,u),(j,\omega)})_{(i,u),(j,\omega)\in
A})\geq bp_{\triangle}(m).$$ Moreover, if $m<{\rm Vol}(\triangle)$
is a turning point of $p_{\triangle}$, and $A\neq
A_m\times\mathbb{Z}/(b)$, then the strict inequality holds.
\end{theorem}
\proof one can show that $\gamma_i=O(\pi^{\lceil\deg(i)\rceil})$.
And, from the equality
$$(\xi_{u}\gamma_{pr-l})^{\sigma^{-1}}=\sum\limits_{\omega=1}^{b}\gamma_{(r,w),(l,u)}\xi_{\omega},$$
we infer that
$$\gamma_{i,j}=O(\pi^{\lceil\deg(pi-j)\rceil}).$$
So we have
$$\sum\limits_{a\in A}{\rm ord}_{\pi}(\gamma_{a,\tau(a)})\geq\sum\limits_{a\in A}\lceil\deg(pa-\tau(a))\rceil\geq
bp_{\triangle}(m).$$ Moreover, if $m$ is a turning point of
$p_{\triangle}$, and $A\neq A_m\times\mathbb{Z}/(b)$, then the
strict inequality holds.\endproof
\section{The Hasse polynomial}
In this section we study $\det(\gamma_{pi-j})_{i,j\in A_m}$.
\begin{definition}For each positive integer $m$, we define $S_m^0$ to be the set of permutations
$\tau$ of $A_m$ satisfying $\tau(0)=0$, and
$$\frac{\tau(a)}{d({\rm sgn}(a))}\geq \deg(pa)-\lceil\deg(pa)-\deg(n)\rceil,\ a\neq0,$$
where $n$ is the element of maximal degree in $A_m\cap({\rm
sgn}(a)\mathbb{N})$.
\end{definition}
\begin{lemma}Let $p>3D$, $m<{\rm Vol}(\triangle)$ a turning point of $p_{\triangle}$, and $\tau$ a permutation of $A_m$.
Then
$$\sum\limits_{a\in A_m}\lceil\deg(pa-\tau(a))\rceil\geq
 p_{\triangle}(m),$$
 with equality holding if and only if $\tau\in S_m^0$.
\end{lemma}
\proof We assume that $M(\triangle) =\mathbb{N}$. The other cases
can be proved similarly. In this case, ${\rm sgn}(a)=+1$ if
$a\neq0$, and the element $n$ of maximal degree in $A_m\cap{\rm
sgn}(a)\mathbb{N}$ is $m-1$. Let $d$ be the nonzero endpoint of
$\triangle$. For $a=0$, we have
$$\lceil\deg(pa-\tau(a))\rceil\geq0$$
with equality holding if and only if $\tau(a)=0$. For $a\neq0$, we
have
$$\lceil\frac{pa-\tau(a)}{d}\rceil\geq\lceil\frac{pa-n}{d}\rceil$$
with equality holding if and only if
$$\frac{\tau(a)}{d}\geq \frac{pa}{d}-\lceil\frac{pa-n}{d}\rceil.
$$
It follows that$$\sum\limits_{a\in
A_m}\lceil\deg(pa-\tau(a))\rceil\geq\sum\limits_{a=1}^n\lceil\deg(pa-n)\rceil=
 p_{\triangle}(m)$$
 with equality holding if and only if $\tau\in S_m^0$.
 The theorem is proved.\endproof

\begin{lemma}We have
$$\gamma_i=\pi^{\lceil\deg(i)\rceil}\sum\limits_{\stackrel{\sum\limits_{j\in\triangle}jn_j=i}
{\sum\limits_{j\in\triangle}n_j =\lceil\deg(i)\rceil}}
\prod_{j\in\triangle}\lambda_{n_j}\hat{a}_j^{n_j}
+O(\pi^{\lceil\deg(i)\rceil+1}).$$  \end{lemma} \proof We have
$$\gamma_i=\sum\limits_{\sum\limits_{j\in\triangle}jn_j=i,n_j\geq0}\pi^{\sum\limits_{j\in\triangle}n_j}
\prod_{j\in\triangle}\lambda_{n_j}\hat{a}_j^{n_j}.$$ We also have
that$$\sum\limits_{j\in\triangle}n_j \geq\lceil\deg(i)\rceil\text{
if }\sum\limits_{j\in\triangle}jn_j=i.$$ The lemma now follows.
\endproof

\begin{definition}For each positive integer $m$, we define
$$H_m(y)=\sum\limits_{\tau\in
S_m^0}\text{sgn}(\tau)\prod_{i\in A_m}
\sum\limits_{\stackrel{\sum\limits_{j\in\triangle}jn_j=pi-\tau(i)}{\sum\limits_{j\in\triangle}n_j
=\lceil\deg(pi-\tau(i))\rceil}}
\prod_{j\in\triangle}\lambda_{n_j}y_j^{n_j}\in\mathbb{Z}_p[y_j\mid
j\in\triangle].$$\end{definition}
\begin{theorem}\label{special-points}Let $p>3D$, and $m<{\rm Vol}(\triangle)$ a
turning point of $p_{\triangle}$. Then
$$\det(\gamma_{pi-j})_{i,j\in A_m}=H_m((\hat{a}_j)_{j\in\triangle})\pi^{p_{\triangle}(m)}+O(\pi^{p_{\triangle}(m)+1/D}).$$
\end{theorem}
\proof Let $S_m$ be the set of permutations of $A_m$. We have
$$\det(\gamma_{pi-j})_{i,j\in A_m}=\sum\limits_{\tau\in S_m}\prod_{i\in A_m}\gamma_{pi-\tau(i)}.$$
The theorem now follows from the last two lemmas.\endproof

\begin{definition}The reduction of $H_m$ modulo $p$ is denoted as
$\overline{H}_m$, and is called the Hasse polynomial of $\triangle$
at $m$.\end{definition}
\begin{theorem}If $p>3D$, and $m<\text{Vol}(\triangle)$ is a turning point of
$p_{\triangle}$, then $\overline{H}_m$ is non-zero.
\end{theorem}
\proof Define $\deg(y_j)=|j|$. Then
$$\prod_{i\in A_m}
\sum\limits_{\stackrel{\sum\limits_{j\in\triangle}jn_j=pi-\tau(i)}{\sum\limits_{j\in\triangle}n_j
=\lceil\deg(pi-\tau(i))\rceil}}
\prod_{j\in\triangle}\lambda_{n_j}y_j^{n_j}$$ has degree
$$\sum\limits_{i\in A_m}|pi-\tau(i)|=\sum\limits_{i\in A_m}{\rm sgn}(i)(pi-\tau(i))$$
$$=p\sum\limits_{i\in A_m}{\rm sgn}(i)i-\sum\limits_{i\in A_m}{\rm sgn}(i)\tau(i)$$
$$\geq p\sum\limits_{i\in A_m}{\rm sgn}(i)i-\sum\limits_{i\in A_m}{\rm sgn}(\tau(i))\tau(i)$$
$$\geq (p-1)\sum\limits_{i\in A_m}{\rm sgn}(i)i,$$
with equality holding if and only if $\tau$ preserves the sign.
Therefore it suffices to show that the reduction of
$$\sum\limits_{\tau\in
S_m^1}\text{sgn}(\tau)\prod_{i\in A_m}
\sum\limits_{\stackrel{\sum\limits_{j\in\triangle}jn_j=pi-\tau(i)}{\sum\limits_{j\in\triangle}n_j
=\lceil\deg(pi-\tau(i))\rceil}}
\prod_{j\in\triangle}\lambda_{n_j}y_j^{n_j},$$ where $S_m^1$
consists of the sign-preserving permutations of $S_m^0$, is nonzero.
One can can prove this by the maximal-monomial-locating technique of
Zhu \cite{Zh1}, as was used by Blache-F\'erard \cite{BF}.\endproof

\begin{definition}We define
$H=\prod_m\overline{H}_m$, where the product is over all turning
points $m<{\rm Vol}(\triangle)$ of $p_{\triangle}$.\end{definition}
\begin{theorem}\label{non-zero}If $p>3D$, then $H$ is non-zero.
\end{theorem}
\proof This follows from the last theorem.\endproof
\section{Proof of the main theorem}
In this section we prove the main theorems of this paper.
\begin{lemma}The Newton polygon of
 $\det(1-\Psi^bs^b\mid
B/\mathbb{Z}_{q}[[\pi^{\frac{1}{D}}]])$ coincides with that of
$\det(1-\Psi s\mid B/\mathbb{Z}_p[[\pi^{\frac{1}{D}}]])$.\end{lemma}
\proof Note that
$$\det(1-\Psi^b s\mid B/\mathbb{Z}_p[[\pi^{\frac{1}{D}}]])
={\rm Norm}(\det(1-\Psi^b s\mid
B/\mathbb{Z}_q[[\pi^{\frac{1}{D}}]])),$$ where Norm is the norm map
from $\mathbb{Z}_q[[\pi^{\frac{1}{D}}]]$ to
$\mathbb{Z}_p[[\pi^{\frac{1}{D}}]]$. The lemma now follows from the
equality
$$\prod\limits_{\zeta^b=1}\det(1-\Psi\zeta s\mid B/\mathbb{Z}_p[[\pi^{\frac{1}{D}}]])
=\det(1-\Psi^b s^b\mid B/\mathbb{Z}_p[[\pi^{\frac{1}{D}}]]).$$
\endproof
\begin{theorem}
The $T$-adic Newton polygon of
 $\det(1-\Psi^bs\mid
B/\mathbb{Z}_{q}[[\pi^{\frac{1}{D}}]])$ is the lower convex closure
of the points
$$(m,\text{ord}_{T}(c_{bm})),\ m=0,1,\cdots.$$\end{theorem}

\proof By the last lemma, the $T$-adic Newton polygon of the power
series $\det(1-\Psi^bs^b\mid B/\mathbb{Z}_{q}[[\pi^{\frac{1}{D}}]])$
is the lower convex closure of the points
$$(i,\text{ord}_{T}(c_i)),\ i=0,1,\cdots.$$
It is clear that $(i,\text{ord}_{T}(c_i))$ is not a vertex of that
polygon if $b\nmid i$. So that Newton polygon is the lower convex
closure of the points
$$(bm,\text{ord}_{T}(c_{bm})),\ m=0,1,\cdots.$$
It follows that the $T$-adic Newton polygon of
 $\det(1-\Psi^bs\mid
B/\mathbb{Z}_{q}[[\pi^{\frac{1}{D}}]])$ is the lower convex closure
of the points
$$(m,\text{ord}_{T}(c_{bm})),\ m=0,1,\cdots.$$
\endproof
 We now prove Theorem
\ref{main1}, which says that $$T-\text{adic NP of }C_{f}(s,T)\geq
{\rm ord}_p(q)p_{\triangle}\text{ if }p>3D.$$

{\it Proof of Theorem \ref{main1}. }Combine the last theorem with
the $T$-adic Dwork's trace formula, we see that the $T$-adic Newton
polygon of
 $C_{f}(s,T)$ is the lower convex closure
of the points
$$(m,\text{ord}_{T}(c_{bm})),\ m=0,1,\cdots.$$
The theorem now follows from the estimate
$$\text{ord}_{\pi}(c_{bm})\geq bp_{\triangle}(m).$$\qed

We now prove Theorem \ref{arithhodge}, which says that
$$p_{\triangle}\geq(p-1)H_{\triangle}^{\infty}$$
with equality holding at the point ${\rm Vol}(\triangle)$.

{\it Proof of Theorem \ref{arithhodge}. }We assume that
$\triangle=[0,d]$. The other cases can be proved similarly. It
suffices to show that
$$p_{\triangle}\geq(p-1)H_{\triangle}^{\infty}\text{ on }[0,d]$$ with equality
holding at the point $d$. Let $0<m\leq d$. We have
$$p_{\triangle}(m)=\sum\limits_{0\leq a<m}(\lceil\frac{(p-1)a}{d}\rceil-\delta_{\in}(a))$$
$$=\sum\limits_{0\leq a<m}(\lceil \frac{pa}{d}\rceil-\lceil \frac{a}{d}\rceil+\delta_{<}(a)-\delta_{\in}(a))$$
$$=\sum\limits_{1\leq a<m}(\lceil \frac{pa}{d}\rceil-1)+\sum\limits_{1\leq a<m}(\delta_{<}-\delta_{\in})(a)$$
$$=\sum\limits_{1\leq a<m}[\frac{pa}{d}]+\sum\limits_{1\leq a<m}(\delta_{<}-\delta_{\in})(a)$$
$$=p\sum\limits_{1\leq a<m}\frac{a}{d}-\sum\limits_{1\leq a<m}\{\frac{pa}{d}\}
+\sum\limits_{\stackrel{1\leq a<m}{d\{\frac{pa}{d}\}\geq m}}1$$
$$=(p-1)H_{\triangle}^{\infty}(m)+\sum\limits_{1\leq a<m}\frac{a}{d}-\sum\limits_{1\leq a<m}\{\frac{pa}{d}\}+\#\{1\leq
a<m\mid d\{\frac{pa}{d}\}\geq m\}.$$ In particular, we have
$$p_{\triangle}(d)=(p-1)H_{\triangle}^{\infty}(d).$$ Note
that
$$\sum\limits_{1\leq a<m}\frac{a}{d}-\sum\limits_{1\leq a<m}\{\frac{pa}{d}\}+\#\{1\leq
a<m\mid d\{\frac{pa}{d}\}\geq m\}$$$$\geq\sum\limits_{1\leq
a<m}\frac{a}{d}-\sum\limits_{\stackrel{1\leq
a<m}{d\{\frac{pa}{d}\}<m}}\{\frac{pa}{d}\}\geq0.$$ It follows that
$$p_{\triangle}(m)\geq(p-1)H_{\triangle}^{\infty}(m).$$ \qed
\begin{theorem}\label{spe} Let $A(s,T)$ be a
$T$-adic entrie series in $s$ with unitary constant term. If
$0\neq|t|_p<1$, then
$$t-adic\text{ NP of
}A(s,t)\geq T-adic\text{ NP of }A(s,T),$$where NP is the short for
Newton polygon. Moreover, the equality holds for one $t$ if and only
if it holds for all $t$.\end{theorem} \proof Write
 $$A(s,T)=\sum\limits_{i=0}^{\infty}a_i(T)s^i.$$
 The inequality follows from the fact that
$a_i(T)\in\mathbb{Z}_q[[T]]$. Moreover, $$t-adic\text{ NP of
}A(s,t)=T-adic\text{ NP of }A(s,T)$$if and only if
$$a_i(T)\in T^e\mathbb{Z}_q[[T]]^{\times}$$ for every turning point
$(i,e)$ of the $T$-adic Newton polygon of $A(s,T)$.  It follows that
the equality holds for one $t$ iff it holds for all $t$.\endproof
\begin{theorem}Let $f(x)=\sum\limits_{u\in\triangle}(a_ux^u,0,0,\cdots)$, and $p>3D$.
If the equality $$\pi_m-\text{adic NP of }C_{f}(s,\pi_m)={\rm
ord}_p(q)p_{\triangle}$$ for one $m\geq1$, then it holds for all
$m\geq1$, and we have $$T-\text{adic NP of }C_{f}(s,T)={\rm
ord}_p(q)p_{\triangle}.$$
\end{theorem}
\proof This follows from Theorems \ref{general-points} and
\ref{spe}.\endproof
\begin{theorem}Let $f(x)=\sum\limits_{u\in\triangle}(a_ux^u,0,0,\cdots)$.
Then  $$\pi_m-\text{adic NP of }C_{f}(s,\pi_m)={\rm
ord}_p(q)p_{\triangle}$$ if and only if
$$\pi_m-\text{adic NP of }L_{f}(s,\pi_m)={\rm
ord}_p(q)p_{\triangle}\text{ on }[0,p^{m-1}{\rm Vol}(\triangle)].$$
\end{theorem}
\proof Assume that $ L_f(s,
\pi_m)=\prod\limits_{i=1}^{p^{m-1}d}(1-\beta_is) $. Then
$$C_f(s,\pi_m )= \prod\limits_{j=0}^{\infty} L_f(q^js, \pi_m)=
\prod\limits_{j=0}^{\infty}\prod\limits_{i=1}^{p^{m-1}d}(1-\beta_iq^js).$$
Therefore the slopes of the  $q$-adic Newton polygon of
$C_{f}(s,\pi_m)$ are the numbers
$$j+{\rm ord}_q(\beta_i),\ 1\leq
i\leq p^{m-1}{\rm Vol}(\triangle), j=0,1,\cdots.$$ One can show that
the slopes of $p_{\triangle}$ are the numbers
 $$j(p-1)+p_{\triangle}(i)-p_{\triangle}(i-1),\ 1\leq i\leq
p^{m-1}{\rm Vol}(\triangle), j=0,1,\cdots.$$ It follows
that
$$\pi_m-\text{adic NP of }C_{f}(s,\pi_m)={\rm
ord}_p(q)p_{\triangle}$$ if and only if
$$\pi_m-\text{adic NP of }L_{f}(s,\pi_m)={\rm
ord}_p(q)p_{\triangle}\text{ on }[0,p^{m-1}{\rm Vol}(\triangle)].$$
\endproof

We now prove Theorems \ref{main2}, \ref{main3} and \ref{expsum}. By
the above theorems, it suffices to prove the following.
\begin{theorem}Let $f(x)=\sum\limits_{u\in\triangle}(a_ux^u,0,0,\cdots)$, and $p>3D$.
Then  $$\pi_1-\text{adic NP of }L_{f}(s,\pi_1)={\rm
ord}_p(q)p_{\triangle}\text{ on }[0,{\rm Vol}(\triangle)]$$ if and
only if $H((a_u)_{u\in\triangle})\neq0$.
\end{theorem}\proof It is known that the
$q$-adic Newton polygon  of $L_{f}(s,\pi_1)$ coincides with
$H_{\triangle}^{\infty}$ at the point ${\rm Vol}(\triangle)$. By
Theorem \ref{arithhodge}, $p_{\triangle}$ coincide with
$(p-1)H_{\triangle}^{\infty}$ at the point ${\rm Vol}(\triangle)$,
It follows that the $\pi_1$-adic Newton polygon of $L_{f}(s,\pi_1)$
coincides with ${\rm ord}_p(q)p_{\triangle}$ at the point ${\rm
Vol}(\triangle)$. Therefore it suffices to show that
$$\pi_1-\text{adic NP of }L_{f}(s,\pi_1)={\rm
ord}_p(q)p_{\triangle}\text{ on }[0,{\rm Vol}(\triangle)-1]$$ if and
only if $H((a_u)_{u\in\triangle})\neq0$.

From the identity $$C_f(s,\pi_1 )= \prod\limits_{j=0}^{\infty}
L_f(q^js, \pi_1),$$ and the fact the $q$-adic orders of the
reciprocal roots of $L_{f}(s,\pi_1)$ are no greater than $1$, we
infer that
$$\pi_1-\text{adic NP of }L_{f}(s,\pi_1) =\pi_1-\text{adic NP of
}C_{f}(s,\pi_1)\text{ on }[0,{\rm Vol}(\triangle)-1].$$Therefore it
suffices to show that $$\pi_1-\text{adic NP of }C_{f}(s,\pi_1)={\rm
ord}_p(q)p_{\triangle}\text{ on }[0,{\rm Vol}(\triangle)-1]$$ if and
only if $H((a_u)_{u\in\triangle})\neq0$. The theorem now follows
from the $T$-adic Dwork trace formula and Theorems
\ref{general-points} and \ref{special-points}.\endproof

\end{document}